\documentclass{article}

\usepackage{a4wide}
\usepackage{amsmath,amssymb}
\usepackage{graphicx,graphics}
\usepackage{subfig}

\newtheorem{definition}{Definition}

\def\a{\alpha}

\begin{document}

\title{Fractional derivatives in Dengue epidemics\thanks{This is a preprint
of a paper accepted for presentation at ICNAAM 2011, 
Numerical Optimization and Applications Symposium,
whose final and definite form will appear in 
\emph{AIP Conference Proceedings}.}}

\author{Shakoor Pooseh\\
Center for Research and Development in Mathematics and Applications\\
Department of Mathematics, University of Aveiro, Portugal\\
\texttt{spooseh@ua.pt}
\and
Helena Sofia Rodrigues\\
School of Business Studies,
Viana do Castelo Polytechnic Institute, Portugal\\
\texttt{sofiarodrigues@esce.ipvc.pt}
\and
Delfim F. M. Torres\\
Center for Research and Development in Mathematics and Applications\\
Department of Mathematics, University of Aveiro, Portugal\\
\texttt{delfim@ua.pt}}

\date{}

\maketitle


\begin{abstract}
We introduce the use of fractional calculus, \textrm{i.e.},
the use of integrals and derivatives of non-integer (arbitrary) order,
in epidemiology. The proposed approach is illustrated
with an outbreak of dengue disease, which is motivated by
the first dengue epidemic ever recorded in the Cape Verde islands
off the coast of west Africa, in 2009. Numerical simulations
show that in some cases the fractional models fit better the reality
when compared with the standard differential models.
The classical results are obtained as particular cases by considering
the order of the derivatives to take an integer value.

\bigskip

\noindent \textbf{Keywords:} fractional calculus, Riemann--Liouville derivatives,
epidemiological model, dengue disease, numerical approximation.

\smallskip

\noindent \textbf{MSC 2010:} 26A33, 34K28, 34K37, 92D30.

\smallskip

\noindent \textbf{PACS:} 02.60.-x, 87.10.Ed
\end{abstract}


\section{Introduction}
\label{introduction}

The incidence of dengue has grown intensely in the last decades and,
according to the World Health Organization, about
40\% of world's population is now at risk. This global
pandemic is attributed to the unprecedented population growth,
the rising level of urbanization without adequate domestic water
supplies, increasing movement of the virus between humans
(due to tourism, migration, or international trade),
and lack of effective mosquito control \cite{Semenza2009}.
Dengue virus is transmitted to humans through the bite of infected
\emph{Aedes} mosquitoes, specially \emph{Aedes Aegypti}. Once
infected, a mosquito remains infected for life, transmitting the
virus to susceptible individuals during feed. Without a vaccine,
vector control remains the only available strategy against dengue.

Appropriate mathematical models can give a deeper insight into
the mechanism of disease transmission \cite{Derouich2006,Otero2008,Thome2010}.
Typically, epidemiologic models are formulated with classical derivatives of integer order.
Here we propose the use of generalized fractional derivatives.
Fractional calculus (calculus of non-integer order),
in spite of its long history as a pure branch of mathematics,
only recently has shown to be useful as a practical tool \cite{MR2736622}.
In this paper we claim that fractional calculus provides
an interesting modeling technique in the context of epidemiology.

We begin by considering a simple epidemiological model
that represents an episode of dengue disease.
The rest of the paper is then dedicated to introduce the notion of fractional
derivative in the sense of Riemann--Liouville,
reformulating the dynamics of the classical model in terms of fractional derivatives,
and finally applying a recent approximate technique
to obtain numerical solutions to the fractional model.
Numerical simulations show that models with fractional derivatives
may provide a better description of reality when compared with the standard ones.


\section{Epidemiological model}
\label{epidemiological model}

In this paper we assume that the host population is divided into
three classes: susceptible, $S_h (t)$, individuals who can
contract the disease; infected, $I_h(t)$, individuals capable of
transmitting the disease to others; and resistant, $R_h (t)$,
individuals who have acquired immunity at time $t$. The total
number of hosts is constant, \textrm{i.e.},
$N_h=S_h(t)+I_h(t)+R_h(t)$. Similarly, we also have two
compartments for the mosquito: $S_m(t)$ and $I_m(t)$ with
$N_m=S_m(t)+I_m(t)$. The model is described
by the system of differential equations
\begin{equation}
\label{ode}
\begin{tabular}{l}
$ \left\{
\begin{array}{l}
\frac{dS_h}{dt}(t) = \mu_h N_h - (B\beta_{mh}\frac{I_m}{N_h}+\mu_h)S_h\\
\frac{dI_h}{dt}(t) = B\beta_{mh}\frac{I_m}{N_h}S_h -(\eta_h+\mu_h) I_h\\
\frac{dR_h}{dt}(t) = \eta_h I_h - \mu_h R_h\\
\frac{dS_m}{dt}(t) = \mu_m N_m -(B \beta_{hm}\frac{I_h}{N_h}+\mu_m) S_m\\
\frac{dI_m}{dt}(t) = B \beta_{hm}\frac{I_h}{N_h}S_m -(\mu_m) I_m
\end{array}
\right. $
\end{tabular}
\end{equation}
subject to given initial conditions
$S_h(0)$, $I_h(0)$, $R_h(0)$, $S_m(0)$ and $I_m(0)$.
The recruitment rate of human and vector populations are denoted as
$\mu_hN_h$ and $\mu_mN_m$, respectively. The natural death rate
for humans and mosquitoes is described by the parameters $\mu_h$
and $\mu_m$, respectively. We assume that $B$ is the average daily
biting (per day) of the mosquito whereas $\beta_{mh}$ and
$\beta_{hm}$ are related to the transmission probability (per
bite) from infected mosquitoes to humans and vice versa. The
recovery rate of the human population is denoted by $\eta_h$.


\section{Fractional Calculus}
\label{fractional calculus}

Fractional calculus was originated at the end of the seventeenth
century and consists in the study of derivatives and integrals
of arbitrarily real or even complex order. Recently,
fractional calculus is experiencing an intensive progress
in both theory and applications \cite{MR2736622,MR2218073}.
Main claim is that a fractional model can give a more realistic
interpretation of natural phenomena. Among several different definitions
that can be found in the literature for fractional derivative,
one of the most popular is the Riemann--Liouville derivative.

\begin{definition}
\label{Def}
Let $x(\cdot)$ be an absolutely continuous function in $[a,b]$ and $0\leq \a<1$.
The (left) Riemann--Liouville fractional derivative of order $\a$,
$_aD_t^\a$, is given by
$$
_aD_t^\a x(t)=\frac{1}{\Gamma(1-\a)}\frac{d}{dt}\int_a^t (t-\tau)^{-\a}
x(\tau)d\tau, \quad t\in [a,b].
$$
\end{definition}
We reformulate system \eqref{ode} using fractional derivatives.
To this end, we simply substitute the first-order derivatives
by Riemann--Liouville derivatives of order $\a$:
\begin{equation}
\label{fde}
\begin{tabular}{l}
$ \left\{
\begin{array}{l}
_0D_t^\a S_h(t) = \mu_h N_h - (B\beta_{mh}\frac{I_m}{N_h}+\mu_h)S_h\\
_0D_t^\a I_h(t) = B\beta_{mh}\frac{I_m}{N_h}S_h -(\eta_h+\mu_h) I_h\\
_0D_t^\a R_h(t) = \eta_h I_h - \mu_h R_h\\
_0D_t^\a S_m(t) = \mu_m N_m -(B \beta_{hm}\frac{I_h}{N_h}+\mu_m) S_m\\
_0D_t^\a I_m(t) = B \beta_{hm}\frac{I_h}{N_h}S_m -(\mu_m) I_m .
\end{array}
\right. $
\end{tabular}
\end{equation}
Note that when $\alpha \rightarrow 1$ the fractional system \eqref{fde}
reduces to \eqref{ode}. Given a certain reality,
our goal is to find the order $\alpha$ that
makes the model more realistic. As we shall see in the next section,
the best value of $\alpha$ in system \eqref{fde}
is in general different from one, \textrm{i.e.},
the classical model \eqref{ode} is often not the best choice.


\section{Numerical simulations}
\label{numerical}

\begin{figure}[!ht]
\begin{minipage}[b]{.5\linewidth}
\centering
\includegraphics[scale=.55]{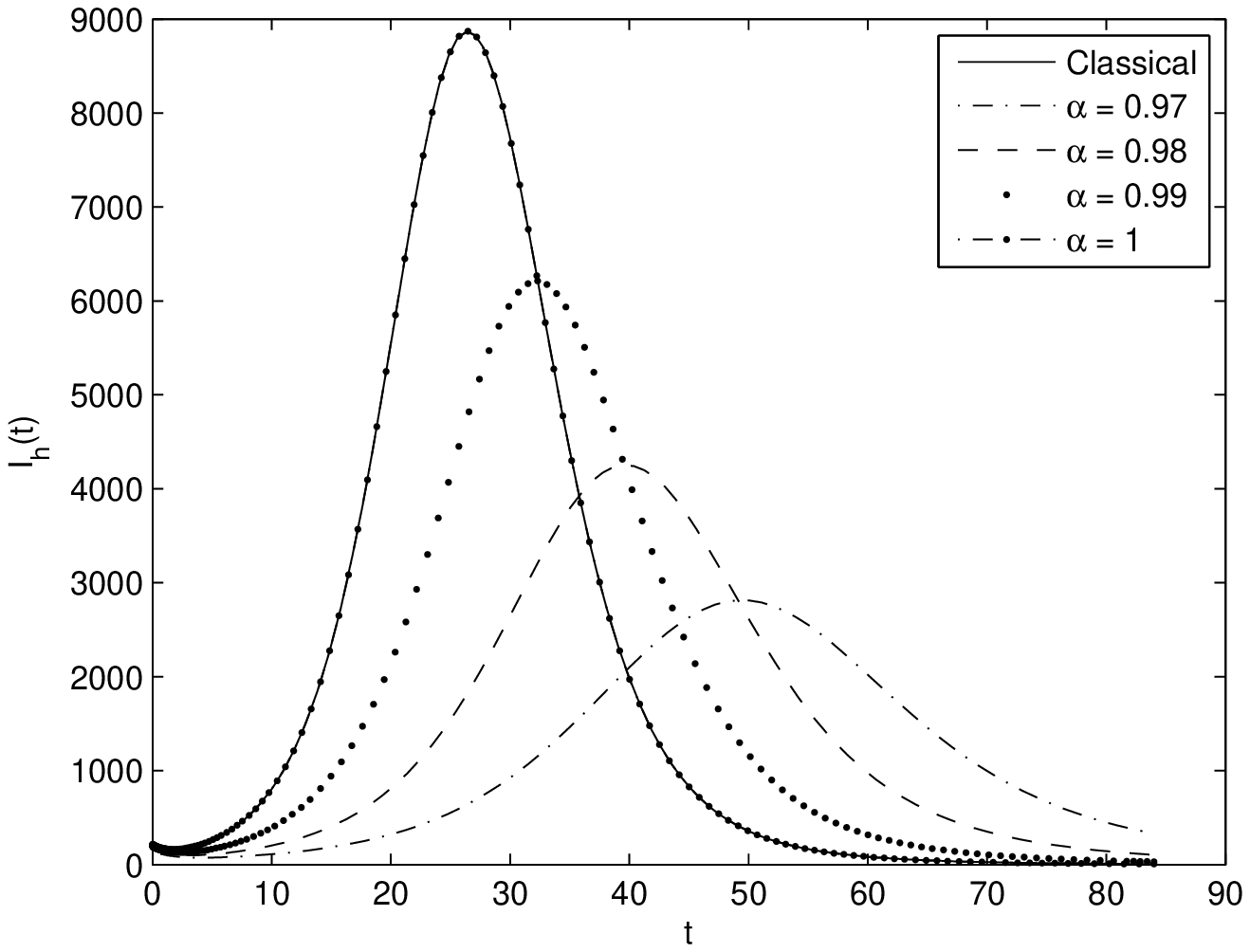}
\caption{solution to the classical model \eqref{ode}
versus the solution to the fractional model \eqref{fde} with
different values of $\a$.} \label{Fig1}
\end{minipage}
\hspace{.3cm}
\begin{minipage}[b]{.5\linewidth}
\centering
\includegraphics[scale=.55]{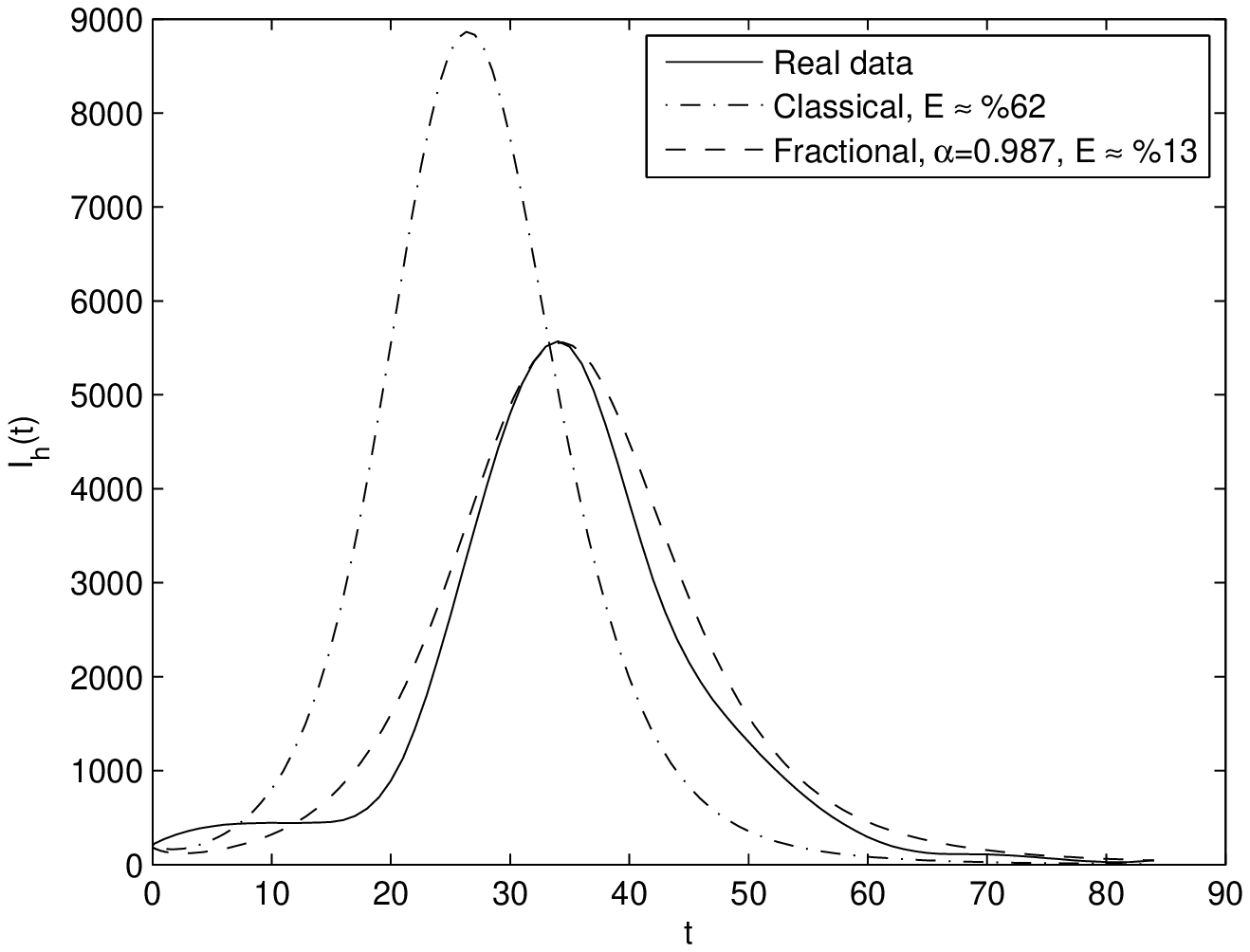}
\caption{real data versus the solution to the classical model \eqref{ode}
and the solution to fractional model \eqref{fde} with $\alpha = 0.987$.} \label{Fig2}
\end{minipage}
\end{figure}

Similarly to the classical theory of differential equations, there
are no general methods to solve systems of fractional differential
equations analytically. Furthermore, the fractional case
is believed to be more difficult to handle even numerically \cite{Kai}.
Recently, an approximation based on a continuous expansion formula for
the Riemann--Liouville fractional derivative has been proposed in
\cite{Atan} and further improved in \cite{APT}. The basic idea is to express
fractional terms by means of a series involving integer-order
derivatives. The immediate result is that one can then use classical methods
to obtain an approximate solution to the original fractional problem.
With the same assumptions of Definition~\ref{Def},
the approximation for $_aD_t^\a x(\cdot)$ is
\begin{equation}
\label{Approximation}
_aD_t^\a x(t)\simeq A(\a,N)(t-a)^{-\a}x(t)+A'(\a,N)
(t-a)^{1-\a}\dot{x}(t)-\sum_{p=2}^N C(\a,p)(t-a)^{1-p-\a}V_p(t),
\end{equation}
where
$$
\left\{
\begin{array}{l}
\dot{V}_p(t)=(1-p)(t-a)^{p-2}x(t)\\
V_p(a)=0, \qquad p=2,3,\ldots,N,
\end{array}
\right.
$$
and $A=A(\a,N)$, $A'=A'(\a,N)$ and $C_p=C(\a,p)$ are given by
\begin{eqnarray*}
A(\a,N)&=&\frac{1}{\Gamma(1-\a)}\left[1+\sum_{p
=2}^N\frac{\Gamma(p-1+\a)}{\Gamma(\a)(p-1)!}\right],\\
A'(\a,N)&=&\frac{1}{\Gamma(2-\a)}\left[1+\sum_{p
=1}^N\frac{\Gamma(p-1+\a)}{\Gamma(\a-1)p!}\right],\\
C(\a,p)&=&\frac{1}{\Gamma(2-\a)\Gamma(\a-1)}
\frac{\Gamma(p-1+\a)}{(p-1)!}.
\end{eqnarray*}
Using \eqref{Approximation}, we approximate system \eqref{fde}
by a system of ordinary differential equations in which,
subject to the order of approximation $N$, new variables of the form
$V_P(\cdot)$ will appear:
\begin{equation}
\label{AppODE}
\begin{tabular}{l}
$ \left\{
\begin{array}{l}
\frac{d\tilde{S}_h}{dt}(t) =
        \left[\mu_h N_h - (B\beta_{mh}\frac{\tilde{I}_m}{N_h}+\mu_h)\tilde{S}_h
        -At^{-\a}\tilde{S}_h+\sum_{p=2}^N C_pt^{1-p-\a}V^{S_h}_p(t)\right]A'^{-1}t^{\a-1}\\
\frac{dV^{S_h}_p}{dt}(t)=(1-p)(t-a)^{p-2}\tilde{S}_h(t), \quad p=2,3,\ldots,N\\
\frac{d\tilde{I}_h}{dt}(t) =
        \left[B\beta_{mh}\frac{\tilde{I}_m}{N_h}\tilde{S}_h -(\eta_h+\mu_h) \tilde{I}_h
        -At^{-\a}\tilde{I}_h+\sum_{p=2}^N C_pt^{1-p-\a}V^{I_h}_p(t)\right]A'^{-1}t^{\a-1}\\
\frac{dV^{I_h}_p}{dt}(t)=(1-p)(t-a)^{p-2}\tilde{I}_h(t), \quad p=2,3,\ldots,N\\
\frac{d\tilde{R}_h}{dt}(t) =
        \left[\eta_h \tilde{I}_h - \mu_h \tilde{R}_h-At^{-\a}\tilde{R}_h
        +\sum_{p=2}^N C_pt^{1-p-\a}V^{R_h}_p(t)\right]A'^{-1}t^{\a-1}\\
\frac{dV^{R_h}_p}{dt}(t)=(1-p)(t-a)^{p-2}\tilde{R}_h(t), \quad p=2,3,\ldots,N\\
\frac{d\tilde{S}_m}{dt}(t) =
        \left[\mu_m N_m -(B \beta_{hm}\frac{\tilde{I}_h}{N_h}+\mu_m) \tilde{S}_m
        -At^{-\a}\tilde{S}_m+\sum_{p=2}^N C_pt^{1-p-\a}V^{S_m}_p(t)\right]A'^{-1}t^{\a-1}\\
\frac{dV^{S_m}_p}{dt}(t)=(1-p)(t-a)^{p-2}\tilde{S}_m(t), \quad p=2,3,\ldots,N\\
\frac{d\tilde{I}_m}{dt}(t) =
        \left[ B \beta_{hm}\frac{\tilde{I}_h}{N_h}\tilde{S}_m -(\mu_m) \tilde{I}_m
        -At^{-\a}\tilde{I}_m+\sum_{p=2}^N C_pt^{1-p-\a}V^{I_m}_p(t)\right]A'^{-1}t^{\a-1}\\
\frac{dV^{I_m}_p}{dt}(t)=(1-p)(t-a)^{p-2}\tilde{I}_m(t), \quad p=2,3,\ldots,N.
\end{array}
\right. $
\end{tabular}
\end{equation}
System \eqref{AppODE} is solved numerically
using the following values that are inspired
in the outbreak occurred in 2009 in Cape Verde: $N_h = 56000$,
$B= 0.7$, $\beta_{mh} = 0.36$, $\beta_{hm} = 0.36$, $\mu_h = 1/(71
\times365)$, $\eta_h = 1/3$, $\mu_m = 1/10$, $m=3$, and
$N_m=m\times N_h$. The initial conditions for the system \eqref{ode}
of ordinary differential equations are: $S_h(0) = Nh-216$,
$I_h(0)=216$, $R_h(0)=0$, $S_m(0)=N_m$ and $I_m(0)=0$. These values,
together with a zero initial condition for $V^{S_h}_p$, $V^{I_h}_p$,
$V^{R_h}_p$, $V^{S_m}_p$, and $V^{I_m}_p$, make system \eqref{AppODE}
to be an ordinary initial value problem that can be treated numerically.
Figure~\ref{Fig1} shows the solution to the systems
\eqref{ode} and \eqref{AppODE} with respect to the variable $I_h$.
The system \eqref{AppODE} is solved with $N=7$ and different values of $\a$.
In order to compare both fractional ($\alpha \in [0,1)$)
and classical ($\alpha = 1$) models,
we include the statistical information from Cape Verde
in Figure~\ref{Fig2}. The best model is, in our case,
the fractional model with $\alpha = 0.987$, which corresponds a
percentage error of thirteen. In contrast, the percentage error
associated with the classical model is sixty two.


\section{Conclusions}
\label{conclusions}

Describing the reality through a mathematical model,
usually a system of differential equations,
is a hard task that has an inherent compromise
between simplicity and accuracy. In this paper we consider
a very basic model to dengue epidemics. It turns out that, in general,
this basic/classical model does not provide enough good results.
In order to have better results, that fit the reality,
more specific and complicated set of differential equations
have been investigated in the literature --- see
\cite{Rodrigues2010a,Rodrigues2010b,Rodrigues2011} and references therein.
Here we propose a completely new approach to the subject.
We keep the simple model and substitute the usual (local)
derivatives by (nonlocal) fractional differentiation. The use
of fractional derivatives allow us to model memory effects,
and result in a more powerful approach to epidemiological models:
one can then design the order $\alpha$ of fractional differentiation
that best corresponds to reality. The classical case is recovered
by taking the limit when $\alpha$ goes to one. Our investigations
show that even a simple fractional model may give surprisingly good results.
However, the transformation of a classical model into
a fractional one makes it very sensitive to the order of
differentiation $\a$: a small change in $\a$
may result in a big change in the final result.

The present work can be extended in several ways:
by fractionalizing more sophisticated models;
by considering different orders of fractional derivatives
for each one of the state variables, \textrm{i.e.},
models of non-commensurate order.


\section*{Acknowledgments}

This work was supported by {\it FEDER} funds through
{\it COMPETE} --- Operational Programme Factors of Competitiveness
(``Programa Operacional Factores de Competitividade'')
and by Portuguese funds through the
{\it Center for Research and Development
in Mathematics and Applications} (University of Aveiro)
and the Portuguese Foundation for Science and Technology
(``FCT --- Funda\c{c}\~{a}o para a Ci\^{e}ncia e a Tecnologia''),
within project PEst-C/MAT/UI4106/2011
with COMPETE number FCOMP-01-0124-FEDER-022690.
Pooseh was also supported by FCT through
the Ph.D. grant SFRH/BD/33761/2009;
Rodrigues by FCT through
the Ph.D. grant SFRH/BD/33384/2008.



\end{document}